\documentclass[%
 reprint,
 amsmath,amssymb,
 aps,
onecolumn,prl]{revtex4}

\usepackage{amsfonts}
\usepackage{amssymb}
\usepackage{graphicx}
\begin{document}

\title{On using Extreme Values to detect global stability thresholds in multi-stable systems:\\
The case of transitional plane Couette flow}

\date{Published in Chaos Solitons and Fractals: http://dx.doi.org/10.1016/j.chaos.2014.01.008}
\author{Davide Faranda}
 \email{davide.faranda@cea.fr}
\affiliation{Laboratoire SPHYNX, Service de Physique de l'Etat Condens\'e, DSM,
CEA Saclay, CNRS URA 2464, 91191 Gif-sur-Yvette, France
}

\author{Valerio Lucarini}
\affiliation{Klimacampus, Universit\"at Hamburg, Grindelberg 5, 20144, Hamburg, Germany. Department of Mathematics and Statistics, University of Reading, UK}

\author{Paul Manneville}
\affiliation{LadHyX, \'Ecole Polytechnique, F-91128 Palaiseau, France}

\author{Jeroen Wouters}
\affiliation{Klimacampus, Universit\"at Hamburg, Grindelberg 5, 20144, Hamburg, Germany\\}

\begin{abstract}
Extreme Value Theory (EVT) is exploited to determine the global stability threshold  $R_{\rm g}$  of plane Couette flow --the flow of a viscous fluid in the space between two parallel plates-- whose laminar or turbulent behavior depends on the Reynolds number $R$.
Even if the existence of a global stability threshold has been detected in simulations and experiments, its numerical value has not been unequivocally defined.
$R_{\rm g}$ is the value such that for $R > R_{\rm g}$, turbulence is sustained, whereas for $R < R_{\rm g}$ it is transient and eventually decays. 
We address the problem of determining $R_{\rm g}$ by using the extremes - maxima and minima -  of the perturbation energy fluctuations.
When $R\gg R_{\rm g}$, both the positive and negative extremes are bounded. 
As the critical Reynolds number is approached from above, the probability of observing a very low minimum increases  causing asymmetries in the distributions of maxima and minima.
On the other hand, the maxima  distribution is unaffected as the fluctuations towards higher values of the perturbation energy remain bounded.
This tipping point can be detected by fitting the data to the Generalized Extreme Value (GEV) distribution and by identifying $R_{\rm g}$ as the value of $R$ such that the shape parameter of the GEV for the minima changes sign from negative to positive.
The results are supported by the analysis of theoretical models which feature a bistable behavior. 
\end{abstract}

\maketitle

\sloppy

\section{Introduction\label{S-Intro}}

The detection of thresholds underlying the sudden shift from one to another dynamical regime, often called {\it critical transitions\/} or {\it tipping points\/}, is a problem of great importance for complex systems, {\it e.g.} global change in climate science, specie extinction in ecology, stresses in materials, etc.
The statistical approach to this question traditionally involve so-called indicators of criticality~\cite{scheffer2009early}.
Some of these indicators are based on modifications of the auto-correlation properties of specific observables when parameters controlling the system approach some critical value, others on the fact that an increase of the variance and the skewness is observed when moving towards tipping points \cite{kuehn2011mathematical}.
Systems alluded to above involve a multitude of ``agents'' acting on widespread spatiotemporal scales.
Understanding their transitions may involve modeling issues which hinder a neat interpretation of the warnings that one could draw from these indicators.
It should therefore be worth studying systems from physics displaying analogous features but where these issues would be kept at a minimal level.
Far-from-equilibrium transitions in macroscopic systems and especially the transition to turbulence in simple flows configurations offer such an opportunity.
With respect to the latter, Navier--Stokes equations are indeed well-known and the flow regime essentially depends on the rate of shear measured by the Reynolds number $R$, with laminar flow for $R\ll1$ and turbulence for $R\gg1$.

The transition to turbulence follows one or the other of two broad routes depending on whether or not the laminar base flow profile displays an inflection point~\cite{Ma10book}.
When an inflection point is present (shear layer, jets, wakes,\dots) the complexity of the fluid increases progressively and there is no marked hysteresis when $R$ is increased or decreased.
 Of particular interest in the present context, the opposite case with no inflection point -- Poiseuille flow in a pipe or a channel, simple shear (Couette) flows between plates or cylinders, boundary layer flow -- is characterized by an abrupt transition with a wide hysteresis range in $R$ called the {\it transitional range\/}.
 It is marked with the coexistence in space of laminar and turbulent domains: turbulent {\it spots\/} in channel and boundary layer flow, turbulent {\it puffs\/} and {\it slugs\/} in pipe flow.
 The lower end of the transitional range is conceptually pinpointed by a {\it global\/} stability threshold $R_{\rm g}$ below which perturbations of arbitrary shapes and amplitudes all decay asymptotically in time.
The upper end of the transitional range is where laminar/turbulent coexistence disappears  and a {\it featureless\/} turbulent regime prevails. This limit may be more fuzzy but a threshold $R_{\rm t}$ can sometimes be identified~\cite{PD05}. 
 
Recently, important progress has been achieved about the determination of $R_{\rm g}$ in Poiseuille pipe flow and the transition threshold attributed to a balance of puff splitting that propagates turbulence by puff decay that eradicates it at this well-defined Reynolds number~\cite{Aetal11}.
The simplicity of the interpretation was mainly due to the quasi-one-dimensional feature of the system along the tube axis and the straightforward interpretation of the stochastic processes involved.
The case of the transition in plane flows is more complicated owing to the quasi-two-dimensional character of the geometry and the greater variety of the local processes entailed in the growth/decay of turbulent domains~\cite{Ma12}.
Accordingly this give some value to all-purpose statistical approaches apt to signal and quantify the proximity of regime changes with large impact, such as the laminar$\,/\,$turbulent transition.

Such studies bear on the time series of given local or global observables.
Below, we restrict ourselves to the consideration of simple shear flow taking place between parallel plates in relative translation.
More information about this system and its transition to turbulence will be given in Section~\ref{S-PCF}.
The distance to laminar flow is then a global variable of particular interest since it clearly discriminates the laminar state where it is identically zero from the turbulent state where it is non-zero by construction.
This observable will be extracted from large scale direct numerical simulations (DNS) of the Navier--Stokes equations in the corresponding geometry (\S\ref{Ss-num}). 
Although the variance and skewness of the times series are  quantities very straightforward to compute, they are of limited value in view of the threshold determination since there is no {\it a priori\/} way to relate their variation to the position of the tipping point.
Complementarily, we will thus develop a statistical approach based on the Extreme Value Theory and propose a criterion allowing the determination of $R_{\rm g}$ with possibly broader applicability.

As will be recalled in Section~\ref{S-EVT}, extremes are distributed according to one out of three possible asymptotic laws forming a single family called Generalized Extreme Value distribution.
The family is parameterized by a {\it shape exponent\/} describing the distribution's tail, either exponential (type 1 or Gumbel law), decaying as a power law (type 2 or Frechet law), or bounded (type 3 or Weibull law).
We will argue that a system approaching a tipping point generically explores the vicinity of a repeller separating the two competing attracting states involved in the bifurcation, so that the statistics of the observable changes from bounded (type 3) away from the tipping point to unbounded/exponentially decaying (type 1) exactly at the tipping point, to unbounded with a fat (power law) tail, {\it i.e.} a Fr\'echet law, beyond the tipping point, even if the transition is not actually observed due to finite observation time.

The main difficulty of the approach lies in uncertainties related to the finiteness of the time series.
This limitation linked to the experimental/computational load at the data production stage cannot be easily overcome.
Our aim is to present a method for the analysis of finite time series rather than an asymptotic theory which, in the case analyzed and in many others of relevant scientific interest, is often inapplicable.
So, in order to support the methodology, we shall rather turn to appropriate toy-models based on theoretical considerations and  designed to reproduce the phenomenology, but for which data collection is no longer a problem and sensitivity studies may be performed.
This will be developed in Section~\ref{S-models}.

Beyond the issues strictly related to the behavior of plane Couette flow, the goal is to provide a new technique for getting reliable information about early warning of critical transitions in complex systems.
The reduction to theoretical models, as a well established way of understanding complex features in a rather simple set up,  might help to cope with technological issues related to producing long data samples.
The example analyzed here may represent a gateway for applications in climate science or ecology. Accordingly, after a summary of our results, we indicate possible research lines for further developments in Section~\ref{S-Conc}.

\section{Overview of Extreme Value Theory\label{S-EVT}}

We start by introducing the asymptotic theory for extremes of statistically independent and identically distributed (i.i.d.) variables, presenting the generalized extreme value distribution and the resulting three possible extreme value laws. Afterwards,  we consider issues arising when the data in the form of finite length time series as the output of some dynamical system.

\subsection{Asymptotic Extreme Value Theory for i.i.d. variables\label{Ss-iid}}

Classical Extreme Value Theory (EVT) states that, under general assumptions, the statistics of maxima $M_m=\max\{ X_0,X_1, ..., X_{m-1}\}$ of i.i.d. variables $X_0, X_1,\dots, X_{m-1}$, with cumulative distribution function  $F(x)$ in the form:
$$
F(x)=P\{a_m(M_m-b_m) \leq x\},
$$
where  $a_m$ and $b_m$ are normalizing sequences, asymptotically obeys  --if there is convergence to a non-degenerate distribution-- a Generalized Extreme Value (GEV) distribution with cumulative distribution function:
\begin{equation}
F_{G}(x; \mu, \sigma,
\kappa)=\exp\left\{-\left[1+{\kappa}\left(\frac{x-\mu}{\sigma}\right)\right]^{-1/{\kappa}}\right\}
\label{cumul}
\end{equation}
with $1+{\kappa}(x-\mu)/\sigma>0 $  \cite{haan2006extreme}.  The {\it location parameter\/} $\mu \in \mathbb{R}$ and  the {\it scale parameter\/} $\sigma>0$ in Equation~\ref{cumul} account for the normalization of the data, avoiding the recourse to scaling constants  $a_m$ and $b_m$ \cite{LLR83}. 

Parameter ${\kappa} \in \mathbb{R}$ in (\ref{cumul}) is the {\it shape parameter\/} also called the {\it tail index\/}. Its sign discriminates the kind of tail decay of the parent distribution:
\begin{itemize}
\item When ${\kappa} = 0$,  the distribution is of Gumbel type (type 1).
It is skewed to the right; the location parameter $\mu$  is equal to the mode  but differs from median and mean.
According to Gnedenko \cite{gnedenko1943distribution}, it is the asymptotic Extreme Value Law (EVL) to be expected when the parent distribution shows an exponentially decaying tail, which includes the normal, log-normal, gamma, or exponential types.
\item The Fr\'echet distribution (type 2), with $\kappa>0$, is observed when the parent distribution possess  a fat tail decaying  as a power law, {\it e.g.} when the bulk statistics obeys a  Cauchy or a T-Student distribution.
\item The Weibull distribution (type 3), with $\kappa<0$, corresponds to a parent distribution having a finite upper endpoint $x_{\rm up} = \mu - \sigma/\kappa$.
The traditional definition of the Weibull distribution relates to minima of variables with a cumulative distribution bounded from below so that the definition above could rather be termed \textit{reversed} Weibull.
\end{itemize}
When properties of maxima and minima are of interest, respectively corresponding to the exploration of the right or left tails of the parent distribution, they can be treated on an equal footing by considering the minima as maxima of the variables after sign reversal \cite{coles2001introduction}. 

\subsection{Asymptotic Extreme Value theory for dynamical systems\label{Ss-DS}}

In the past decade significant progresses have been made in understanding EVLs related to the output of a dynamical systems.
Generally speaking this output is in the form of time series of some observable computed from orbits  \cite{balakrishnan1995extreme} and difficulties arise from translating results for i.i.d. variables to such an observable evaluated all along a typical trajectory that  samples the natural measure, while being autocorrelated {\it via\/} the underlying dynamics.
As shown by Leadbetter {\it et al.} \cite{LLR83}, the independence condition can be relaxed on general grounds and replaced by appropriate {\it mixing conditions\/}  expressing the decay of correlations.
In the case of deterministic dynamical systems, such a mixing is naturally expected to stem from the sensitivity to initial conditions inherent in chaos and convergence to EVLs was indeed obtained for a special class of  observables by Freitas {\it et al.\/} \cite{FFT10}.

Another issue specific to dynamical systems, deterministic or noisy, is related to the clustering of extremes, as discussed originally by Newell \cite{newell1964asymptotic}, Loynes \cite{loynes1965extreme}, or O'Brien \cite{o1974limit}, and thoroughly explained in \cite[\S3.7]{LLR83}.
In such a situation, introducing the {\it extremal index\/} $\theta$, convergence to EVLs is obtained for $F_{G}(x; \mu, \sigma,\kappa)^\theta$, where $F_{G}$ is the GEV distribution introduced earlier, whereas  $1/\theta$ is a measure of the cluster size.
Alternatively, theorem 3.7.2 in \cite{LLR83} can be exploited to hide the dependence on $\theta$ by appropriately adjusting $\mu$ and $\sigma$, which will be done automatically when fitting $F_{G}$ against the empirical data from the numerical simulations.

Recent advances related to the use of Extreme Value Theory for dynamical systems have been mainly limited to observables  sampling a typical trajectory while making reference to specific a point on the attractor \cite{FFT10,FFT11, lucarini2012extreme, faranda2013recurrence}.
This kind of observable is however experimentally out of reach when dealing with systems with many degrees of freedom.
For such systems, especially those governed by partial differential equations that cannot be brought to (very) low dimension by inertial manifold reduction~\cite{Te90}, the computation of trajectories is indeed computationally too expensive and the existence of the attractor mainly a conceptual view.
The behavior of the system is then often depicted using {\it physical\/}  (macroscopic) observables, {\it e.g.} total energy or angular momentum,  temperature, {\it etc.}, extreme values of which are analyzed in a pure statistical way, and not from an investigation of its dynamical properties at a local (microscopic) level.

\subsection{Selection procedures and finite samples}
Gnedenko's results relate to selecting extremes using the {\it block maxima approach\/} which consists in dividing the time series  of the considered observable of length $s$ into $n$ bins each containing the same number $m$ of observations, next selecting the maximum (or the minimum) value in each of them \cite{coles2001introduction}.
Then, the GEV distribution is fitted and convergence towards  a member of the GEV family proven in the limits $m\to\infty$, and $n\to\infty$.
When moving from theory to practice, it remains first to show that this asymptotic limit can be reached even with a finite amount of data $s=n\times m$.
Comparing analytical results and empirical findings from numerical experiments, the authors of  \cite{faranda2011numerical}, among whom two of us (D.F. \& V.L.), studied how a robust estimation of parameters could be obtained.
For typical chaotic maps, good agreement between theoretical and experimental parameters was obtained for $n\gtrsim1000$ and $m\gtrsim1000$, which will here help us to fix the order-of-magnitude of the amount of data needed.

\subsection{Extremes value laws near crisis points}

Let us consider a dynamical system, either deterministic or noisy, controlled by some parameter  $\lambda$ which,
when increased beyond some value  $\lambda_{\rm crit}$, drives it through a {\it critical transition\/}, making it tumble from one operating point to another.
Here the word `critical' has the meaning it takes, say, in environmental sciences, where the expression {\it tipping point\/} is also used.
In dynamical systems theory, one would speak of a {\it saddle-node\/} bifurcation or some appropriate generalization of it, namely {\it crises\/}~\cite{Ot93}.%
\footnote{In statistical physics, the present situation would correspond to a (discontinuous) {\it first-order\/} phase transition, {\it e.g.} a liquid-gas transition, as opposed to a (continuous) {\it second-order\/} phase transition, {\it e.g.} a para-ferromagnetic transition, studied within the framework of {\it critical phenomena\/}, where the word `critical' thus gets a different meaning  through the definition of {\it critical exponents\/} and related {\it universality classes\/}~\cite{St88}.\label{FN1}}

On general grounds we may expect that physical observables have bounded fluctuations and that their extremes follow Weibull distributions~\cite{holland2012extreme,lucarini2014towards}.
Gaussian fluctuations (featuring Brownian motion of microscopic degrees of freedom) would yield the formal possibility of infinite extremes and thus Gumbel distributions, but the convergence towards this law is logarithmically slow \cite{hall1979rate}, which makes it unobservable in practice.

This viewpoint is however relevant only as long as there is only one typical time scale in the system and the asymptotic time behavior corresponds to a single component attracting set on this time scale, which is explicitly stated in the formulation of mathematical theorems when the existence of an invariant probability measure is assumed, but most of the time tacitly admitted.
If two attracting pieces are in competition and that, either under the effect of external noise or due to internal chaotic fluctuations, two time scale are present, a short one related to transitive dynamics within an attracting component and a long one corresponding to intermittent jumps from one to the other component, the picture has to be modified.

On the short time scale, we have the situation considered up to now of a Weibull distrubution. 
Noise, either extrinsic or intrinsic, does not change the picture:
It adds the possibility of extremes that have to be considered as {\it local\/} at the scale of that part of the attractor which is visited with probability one for a time series of the typical duration $s$ used in the evaluation of the shape of the distribution, which makes the Gumbel distribution irrelevant in practice due to logarithmic convergence.

On the long time scale,  some extremes correspond to noisy excursions directed toward the saddle-state and gain a {\it global\/} status as they can trigger jumps from one to the other component.
The probability increases that the observable visits corresponding ``anomalous'' values associated to these global extremes during a time series of length $s$.
With the tail of the parent distribution becoming heavier and heavier, the EVL will turn to a Fr\'echet distribution.
Upon continuous variation of the control parameter, the shape exponent can thus be expected to cross zero from negative values, and $\lambda_{\rm crit}$ be defined as the value at which this crossing happens.  

In high dimensional systems of interest, critical transitions are featured by specific physical observables that experience abrupt changes when the system crosses its tipping point.
Like in the above discussion more directed to low-dimensional systems, these observables will display deviations of greater amplitude in the direction of the state the system is doomed to tumble, than in the opposite direction.
This generally implies an increase in the skewness of the distribution, backing the suggestion that this quantity could be an early warning indicator of  a tipping point~\cite{guttal2008changing}.
The skewness being a measure of the asymmetry of a distribution, the method performs well when the distribution is symmetric for $\lambda\ll\lambda_{\rm crit}$ but may fail if the distribution is already badly skewed.
The extreme value analysis of maxima and minima of the time series is able to overcome this problem and locate the critical transition.

It is generally understood that, even if no real ``energy'' can be defined in a strict sense, the system evolves in some {\it energy landscape\/}.
As long as the considered operating point is stable against any perturbation, the geography of this landscape is simple, with a single minimum (sink) and trajectories bouncing around it.
When several different operating points coexist, each being locally stable and having its own attraction basin, trajectories gain some probability to crawl over some pass in that landscape from one sink to another.
Variation of the control parameter is then naively understood to produce changes of relative  altitudes of sinks, bumps and saddles.
A critical transition taking place at some unknown $\lambda_{\rm c}$, as long as $\lambda$ is far below $\lambda_{\rm c}$ the probability of tipping during an experiment producing a time series of length $s$ is completely negligible and the system is assigned to a given sink that it explores with essentially bounded variations of the observable, leading one to expect some finite and negative shape exponent.
When $\lambda$ increases, the probability for the system to leave its current operating point within time $s$ also increases, supposedly changing the tail of the parent distribution. 
The goal is to locate $\lambda_{\rm c}$ with the same criterion as before bearing on the shape parameter using information extracted from the data series at our disposal.
As can be inferred from the discussion above, the main difficulty lies in adjusting the width $m$ of the bins used in the bock maxima approach, with $m$ sufficiently large to sample anomalous extremes typical of the dangerous excursions toward tipping when keeping $n=s/m$ sufficiently large for reliable statistics, while the length $s$ of the time-series is subjected to the experimental constraints.

In the next section we apply the method to plane Couette flow as a prototype of high-dimensional system accurately described by a set of partial differential equations controlled by a single parameter, but working in a complex regime.

\section{Results for plane Couette flow \label{S-PCF}}
\subsection{The transition to/from turbulence in plane Couette flow}

Plane Couette flow is the prototype of plane flows with laminar velocity profiles deprived of inflection points and transiting to turbulence in a subcritical fashion.
This flow configuration refers to the shearing of a viscous fluid in the space between two parallel plates in relative motion.
The plates, at a distance $2h$, translate in opposite directions at a speed $U$ and the flow results from the viscous drag acting on the fluid with kinematic viscosity $\nu$.
The nature of the flow regime, either laminar or turbulent, is controlled by a single parameter,  the Reynolds number $R=Uh/\nu$.
The laminar flow depends linearly on the coordinate normal to the plates and is known to remain stable against infinitesimal perturbations for all values of $R$, while turbulent flow is instead observed under usual conditions when $R$ is sufficiently large, typically of order 400--500, when increasing $R$ without particular care.

As $R$ is decreased from high values for which the flow is turbulent, a particular regime appears at about $R_{\rm t}\approx410$ where turbulence intensity is modulated in space~\cite{PD05}.
When the experimental setup is sufficiently wide, a pattern made of oblique bands, alternatively laminar and turbulent, becomes conspicuous.
Bands have a pretty well defined wavelength and make a specific angle with the streamwise direction.
As $R$ is further decreased, they break down and leave room to the laminar base flow below $R_{\rm g}\approx 325$. 
Experiments show that the streamwise period%
\footnote{By convention, the streamwise direction is  along $x$, the normal to the moving plates defines the $y$ direction, and $z$ denotes the spanwise direction.\label{FN2}}
$\lambda_x$ of the band pattern is roughly constant ($\lambda_x\simeq110h$) while the spawise period $\lambda_z$ increases from about $55h$ close to $R_{\rm t}$ to about $85h$ as $R$ decreases and approaches $R_{\rm g}$~\cite{PD05}.
Whereas the turbulence self-sustainment process in wall-bounded flows is well understood~\cite{Wa97}, the mechanisms explaining band formation are still somewhat mysterious.

The transition thus display a large amount of hysteresis. A similar situation is to be found in several other flow configurations, circular Couette flow, the Couette flow sheared by coaxial cylinders rotating in opposite directions, plane channel, the flow between two plates driven by a pressure gradient, as well as in Poiseuille flow in a circular tube, see Section~\ref{S-Intro}.

\subsection{Conditions of the numerical experiment\label{Ss-num}}

The transition to turbulence in plane Couette flow has been studied numerically by a number of authors.
System sizes required to observe the oblique band regime in  Navier--Stokes DNSs are numerically  quite demanding~\cite{Detal10}.
In order to  reduce the computational load, Barkley \& Tuckerman performed their computations in a cleverly chosen narrow but inclined domain~\cite{TB11}.
The drawback is however to freeze the orientation beforehand, forbidding any angle or orientation fluctuation.
Previous work by one of us (P.M.) has shown that another way to decrease computer requirements was to accept some under-resolution of the space dependence, especially in the wall-normal direction $y$~ \cite{MR11}.
All qualitative features of the transitional range are indeed well reproduced in such a procedure, including orientations fluctuations.
Quantitatively, the price to pay appears to be a systematic downward shift of the $[R_{\rm g},R_{\rm t}]$ interval (see below), giving supplementary evidence of the particular robustness of the band regime.

In this work, DNSs have been performed in a domain of constant size able to contain one pattern wavelength in each direction, i.e. $(L_x,L_z)\equiv(\lambda_x\times\lambda_z)$, with $L_x=108$ and $L_z=64$.
This size seems well adapted to the central part of the transitional domain, {\it i.e.} slightly too wide for $R\approx R_{\rm t}$ and slightly too narrow at $R\approx R_{\rm g}$, with mild consequence on the effective value of these thresholds, as guessed from a Ginzburg--Landau approach to this pattern forming problem~\cite{RM11a}.
Such finite-size effects~\cite{PM11} also account for the intermittent reentrance of featureless turbulence.%
\footnote{Intermittent reentrance of featureless turbulence was present in the Barkley--Tuckerman simulations for identical reasons.}

The well-validated open-source software \textsc{ChannelFlow} \cite{gibson2007channelflow} has been used throughout the study.
This Fourier--Chebyshev--Fourier pseudo-spectral code is dedicated to the numerical simulation of flow between parallel plates with periodic in-plane boundary conditions.
In the wall-normal direction (see Note~\ref{FN2}), the spatial resolution is a function of the number $N_y$ of Chebyshev polynomials used.
The in-plane resolution depends on the numbers $(N_x, N_z)$ of collocation points used in the evaluation of the nonlinear terms.
From the 3/2 rule applied to remove aliasing, this corresponds to solutions evaluated in Fourier space using $\frac23 N_{x,z}$ modes, or equivalently to effective space steps $\delta^{\rm eff}_{x,z}=\frac32 L_{x,z}/N_{x,z}$.
Numerical computations have been performed using three different resolutions: low ($N_y=15$, $N_x=L_x$, $N_z=3L_z$), medium ($N_y=21$, $N_x=2L_x$, $N_z=6L_z$), and high ($N_y=27$, $N_x=3L_x$, $N_z=6L_z$) for which we expect $[R_{\rm g},R_{\rm t}]\approx[275,350]$, $[300,380]$, and $[325,405]$, respectively; see Fig.~6 in  \cite{MR11}. 
Simulations have been performed chronologically from low to high resolution, confirming these transitional range estimates but results will be presented in the reverse order since they require less and less computing power, which allows better and better statistics at comparable numerical load.

\subsection{The EVT analysis of turbulent energy near the band-breaking point}

In this section we show results about the changes in the  extreme value distributions of quantity $E_{\rm t}$ defined as the mean-square of the perturbation velocity $\mathbf{\tilde v}$, the difference  between the full velocity field $\mathbf{v}$ and the base flow velocity $\mathbf{v}_{\rm b} = y\, \mathbf{e}_x$.
Physically speaking, apart from a factor $\frac12$, it is thus the kinetic energy contained in the perturbation and accordingly a good measure of the distance to laminar flow where it  is identically zero.
Here we focus on the determination of $R_{\rm g}$ using extremes as sketched above.
A similar approach could be developed to study the behavior of extremes associated to turbulence reentrance around $R_{\rm t}$ but the subcritical character of the transition at this threshold value is still unclear.

For each value of the Reynolds number, very long simulations are performed and, once the time series of $E_{\rm t}$ has reached a stationary state, maxima (minima) are extracted in bins of fixed block length as described in the previous section.
We then fit the maxima (minima after sign change) to the GEV distribution by using a Maximum Likelihood estimation as described in \cite{faranda2011numerical}.
The results could have equivalently been obtained by using other estimators.
The choice of the bin length $m$ is crucial: in the asymptotic regime the value of  the shape parameter should be independent of $m$.
We have tested that, within the confidence intervals, this happens for $m>1000$.

\begin{figure}
\begin{center}
\includegraphics[width=0.8\textwidth]{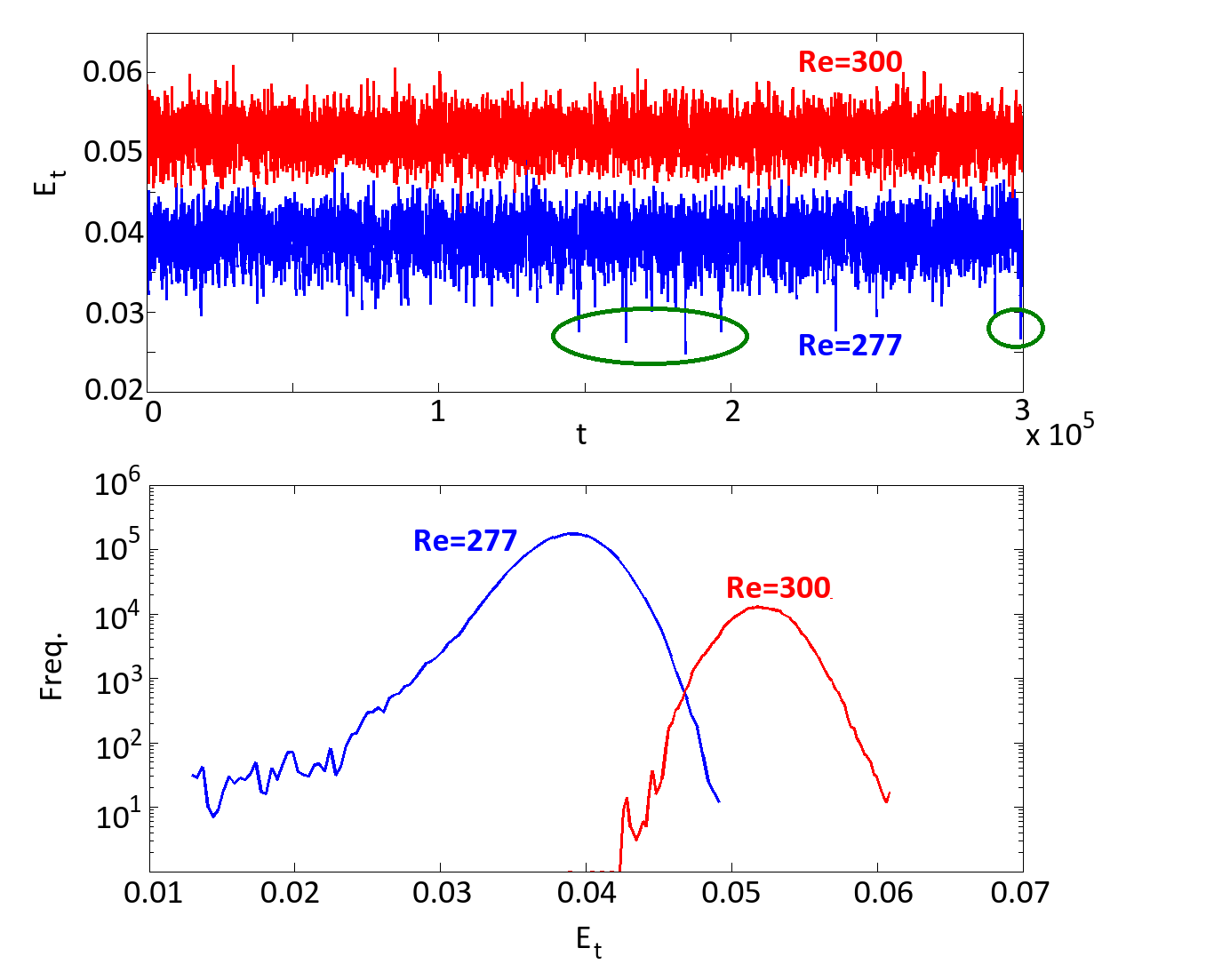}
\caption{(Color online) Perturbation Energy $E_{\rm t}$ for two simulations at low resolution.
Upper panel:  $E_{\rm t}$ as a function of time.
Lower panel: Histograms of $E_{\rm t}$ in log-linear scale. \label{couette}}
\end{center}
\end{figure}

The shape parameter is next analyzed as a function of the Reynolds number.
A first intuition on how the method should work comes from looking at the data series and the histograms shown in figure~\ref{couette}, see caption for details.
The series in red refers to a value of Reynolds inside the band regime ($R=300$),
with fluctuations exploring a limited interval.
The series in blue, with $R$ is fixed just above $R_{\rm g}$ ($R=277$), illustrates a clear tendency to intermittently visit states with very low values of the energy. These events, spotted in the green ovals, crucially contribute to a shift towards Fr\'echet laws since the fit to the GEV returns a Weibull EVL when removing them from the histogram. 

\begin{figure}
\begin{center}
\includegraphics[width=1.0\textwidth]{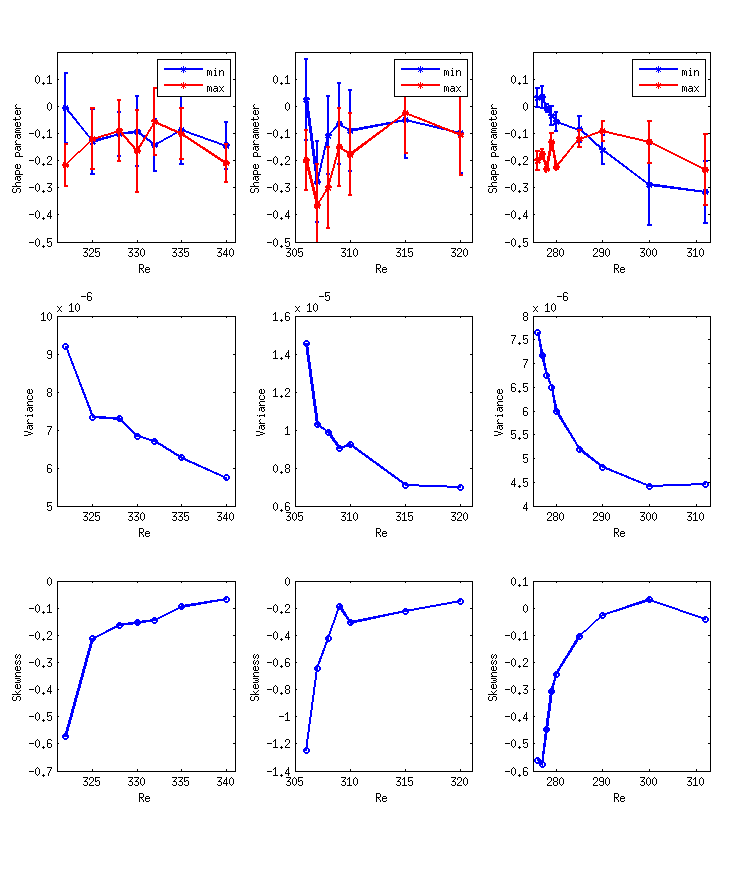}
\caption{Tipping-point indicators for plane Couette flow as functions of $R$.
Upper panel: Shape parameter $\kappa$; red: maxima, blue: minima; error bars represent 95\% confidence intervals, $m=1000$. Center and bottom panels: Variance and skewness of the full series, respectively. Left: High resolution. Center: Medium resolution. Right: Low resolution.
\label{couetteTR}}
\end{center}
\end{figure}

Let us start with the localization of the global stability threshold $R_{\rm g}$ in simulations performed at high resolution, namely  $N_x=216, N_y=27, N_z=384$.
Results are shown in Fig.~\ref{couetteTR} (left column) for the shape parameter (upper panel), to be  compared to the two common early warnings indicators based on the bulk statistics: the skewness (middle panel) and the variance (lower panel).
When approaching $R=322$  the shape parameter for the distribution of minima changes its sign, whereas for the maxima it remains negative in agreement with what was stated in the previous sections.
It is however evident that these results need confirmation since a limited set of Reynolds numbers has been studied and a single slightly positive value of $\kappa$ has been obtained for $R=322$, with error bars so large that  the significance of the result is rather limited.
The variance and the skewness of the time series follow what is expected from the statistics of global observables at a tipping point, namely a monotonic trend towards larger values.
Since $E_{\rm t}$ visits lower energy states, the skewness becomes more negative so that only the distribution of minima is affected.
However, as noticed previously, no definite threshold value $R_{\rm g}$ can be inferred from the consideration of the variance and skewness curves. 

 In fact, at full resolution, the cumulated amount of CPU time required to produce series of length $s=2.5\cdot 10^5$  time units  was beyond $10^5$~CPU hours, making it practically impossible to obtain much longer series with the available resources. 
In order to support the so-obtained results, we have exploited the fact that downgrading the resolution preserves the qualitative features of the transition, up to a shift of transitional range~\cite{MR11}.
At  medium resolution $(N_x=216, N_y=21, N_z=384)$, all time series have been stopped at $s=2\cdot 10^5$ time units. In  these conditions, band breakdown was never observed for $R>306$.
The results  shown in Fig. \ref{couetteTR} (center column) confirm those at high resolution, with a slightly more pronounced change of sign of the shape parameter at $R=306$ but point out the need of more and much longer series around the  global stability threshold.

Further downgrading the resolution to  $N_x=108, N_y=15, N_z=192$, allowed us to produce series lasting nearly one order of magnitude longer than above, up to $2\cdot10^6$ time units.
As a matter of fact, by collecting a greater statistics of maxima, the uncertainty on the estimation of the shape parameter could be greatly reduced, as shown in Fig. \ref{couetteTR} (right column).
In view of our proposal to define $R_{\rm g}$ using extreme value statistics, the results at low resolution look much more convincing  then those produced at higher resolutions since a clear monotonic variation of the shape parameter for minima is now observed upon decreasing $R$. As soon as $R_{\rm g}\leq 278$, a Fr\'echet distribution is found for the minima of the turbulent energy.

Values of $R_{\rm g}$ determined here for the different resolutions studied are not much different from those given in \cite{MR11} obtained by inspection of individual cases without any systematic criterion and using much shorter time series.
While the proposed new methodology to define the global stability threshold seems appropriate, it would be interesting to determine a rescaling procedure defining a master curve common to our three cases independently of the resolution and, next, to justify  $\kappa=0$ as signature of the threshold by some heuristic argument.
The complexity of plane Couette flow forbids us to scrutinize these issues by means of DNSs and suggests to make use of simplified toy models, as considered in the next section.

\section{Theoretical models of critical transitions\label{S-models}} 

A good candidate for testing the identification of the global stability threshold using methods based on GEV parameters  is a slightly modified version the model originally introduced in \cite{dauchot1997local}:
\begin{equation}
{\rm d}X/{\rm d}t=-(\mu+u\xi(t))X+Y^2, \qquad {\rm d}Y/{\rm d}t=-\nu Y +X -XY.
\label{bristol1}
\end{equation} 
Here $X$ and $Y$ may be related to the amplitudes involved in the self-sustaining process of turbulence.
Parameters $\mu$ and $\nu$ are damping coefficients accounting for viscous effects and assumed to vary as $1/R$.
Non-linearities preserve the energy $E=\frac12(X^2+Y^2)$ in the same way as the advection term of the Navier--Stokes equations.
Noise is here introduced in a multiplicative way {\it via\/} the term $u\xi(t)$, where $\xi(t)$ is a white noise and~$u$ its amplitude,  as proposed by Barkley~\cite{Barkley2011modeling}.
A saddle-node bifurcation takes place at $\mu\nu=\frac14$. The trivial solution $X=\nobreak Y=\nobreak0$, corresponding to laminar flow, competes with two nontrivial solutions on the interval $\mu\nu=\left[0,\frac14\right]$, the stable nontrivial solution being assimilated to turbulent flow.
Unlike the additive noise considered in~\cite{manneville2005modeling}, the multiplicative noise taken here does not affect the trivial state and can be understood as a fluctuating turbulent-like contribution to effective viscous effects. 
Whenever the system undergoes a transition towards the laminar state, the simulation is restarted from the stable nontrivial fixed point.

Whereas for plane Couette flow only one simulation could be performed at each Reynolds number, here we can easily produce ensembles of realizations for given parameter sets $(\nu, \mu, u)$, extract corresponding GEV shape parameters, and average them over the realizations, from time series of the energy $E$ that, measuring of the distance to the trivial state, remains an appropriate observable.
In view of locating the transition from a probabilistic viewpoint  we now want to relate the conditions when the probability of transition becomes significant to the change of sign of the GEV parameter $\kappa$.

\begin{figure}
\includegraphics[width=1.0\textwidth]{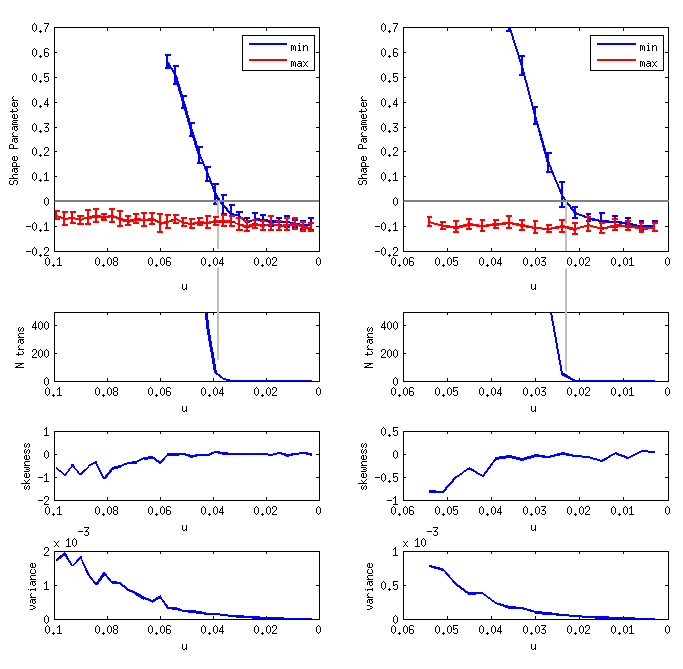}
\caption{Variation of different indicators of critical transitions as functions of the noise intensity for model (\ref{bristol1}) with $\mu=1$.
Left: $\nu=0.2475$. Right: $\nu=0.2487$.
Top row: Averaged shape parameter $\kappa$; red: maxima, blue: minima; error bars represent  the standard deviation over the ensemble of 30 independent realizations.
Second row: Number of transitions observed (see text).
Vertical lines between the two top panels point to the critical value $u_{\rm c}$ for which $\kappa=0$.
Third and fourth row, respectively,  averaged variance and skewness.
Note that $u$ increases to the left and not to the right as usual.}
\label{Bristolm}
\end{figure}

Results of two different simulations are shown in Figure~\ref{Bristolm}.
Here, the control parameter is the intensity of the noise $u$ whereas $\mu$ and~$\nu$ are kept fixed.
The left plots refer to the case ${\rm d}t=0.01$, $\mu=1$, $\nu=0.2487$, $n=10^3$, $m=10^6$, whereas the  right  ones refer to $\nu=0.2475$ with the other parameters left unchanged.
For each values of $u$, ensembles of 30 realizations has been prepared.
The upper panels show the variation with $u$ of the shape parameter averaged over the realizations with error bars corresponding to the standard deviation over each ensemble.
The plots in the second row display the number $N_{\rm tr}$ of times the system has undergone a critical transition during a whole simulation performed at given $u$, {\it  i.e.} the duration of each realizations being $s$, the number $N_{\rm tr}$ is obtained by counting how many times the system has reached the laminar state and been reset to the nontrivial fixed point in $30\times s$ time units.

In both simulations, the shape parameters vary like in plane Couette flow.
For the distribution of minima, it crosses zero and changes its sign when $u$ reaches a certain value $u_{\rm c}$ (vertical lines between the first and second rows in Fig.~\ref{Bristolm}) that depends on the other parameters $\mu$ and~$\nu$, while it remains negative for the maxima.
That the value so assigned to $u_{\rm c}$ acts as an effective transition warning stems from the comparison of the two upper rows  of Fig. \ref{Bristolm}, from which we infer that  $N_{\rm tr}$ markedly increases when $\kappa$ goes through zero.
Experiments have been repeated for other values of $\nu$ (not shown here) with identical results.
The plots in the third  and fourth row respectively refer to the variance and the skewness of the bulk statistics of $E$.
In contrast, these indicators based on the bulk statistics, while displaying the expected trends, do not show any accident allowing us to locate any threshold.

A problem that cannot be explored with sufficient data from the DNS of plane Couette flow relates to whether  GEV distributions are properly fitted, which will be solved if we show that scaling laws relate the bin length $m$ to other parameters of the system while preserving the graph of the shape parameter.
This would indeed imply that, no matter the total length of the time series, the same warning conditions would be obtained provided the data samples provided they are not too short.
This property would be highly desirable in view of  the extension of the present approach to other applications, {\it e.g.} in the environment.
We scrutinize this problem by considering the simplest model  featuring bi-stability in the presence of random noise:
\begin{equation}
{\rm d}X=-V'(X){\rm d}t +\epsilon {\rm d}W
\label{dw}
\end{equation}
with potential  $V(X)=\frac{1}{4}X^4 -aX^2+\lambda X$, where $a>0$, $\lambda>0$  are parameters and  $W$ is a Wiener process of amplitude $\epsilon>0$.
Like parameter $\mu$ in Model (\ref{bristol1}), $\lambda$ serves as a control parameter and pursuing the analogy to plane Couette flow, we can assume $\lambda \propto 1/R$.

We consider System (\ref{dw}) for values of $\lambda$ such that, in the deterministic limit, it features two stable fixed points ($\bar X_1<0$ and $\bar X_2>0$) and an unstable fixed point $\tilde X$.
The asymptotic behavior of this system can be assessed in terms of the solution to a Fokker--Plank equation~\cite{risken1989fokker}.
Here we are rather only interested in the finite-time behavior and we restrict to simulations such that the noise is sufficiently small to leave  the system confined in one of its two wells and does not push it over the saddle. 

\begin{figure}
\begin{center}
\includegraphics[width=0.8\textwidth]{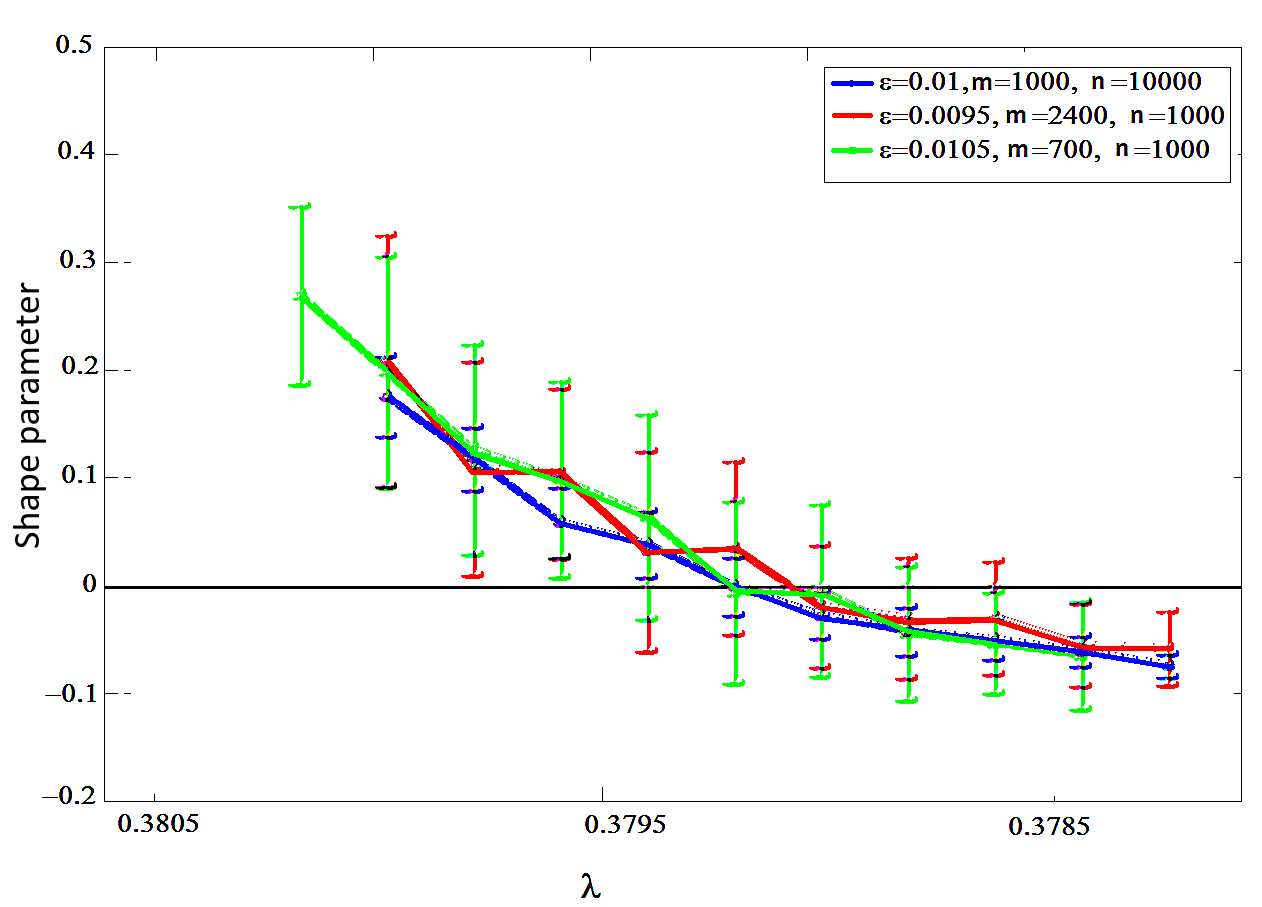}
\caption{Shape parameter $\kappa$  for the minima of Observable $X$ as a function of $\lambda$ for Model~(\ref{dw});   error bars represent the standard deviation over the ensemble.}
 \label{langevin}
 \end{center}
 \end{figure}

In the block maxima approach, extremes are detected in bins of length $m$ whereas $n$ is kept fixed and supposed to be sufficiently large that the resulting GEV distribution is appropriately sampled.
The bin length $m$ is thus the characteristic time to be used.
On the other hand, in System (\ref{dw}), escape from a well takes an average time  $\langle\tau\rangle$  given by
\begin{equation}
\langle \tau \rangle \propto \exp\left(2\Delta V / \epsilon^2\right),
\label{conj}
\end{equation}
where $\Delta V = V\big(\tilde X\big) - V\left(\bar X\right)$, where $\bar X$ is the coordinate of the minimum of $V$ of interest and $\tilde X$ that of the saddle \cite{hanggi1990reaction}. 
Following the same procedure  as  for Model (\ref{bristol1}), we consider the minimum with positive abscissa $\bar X_2>0$ as initial condition and prepare ensembles of 30 realizations for each value of $\lambda$.
The observable chosen here is just variable $X$ itself.
For the minima, the shape parameter $\kappa$ averaged over the realizations is seen to change its sign, as expected since minima correspond to excursions toward the basin of attraction of the other minimum, {\it i.e.} aborted transitions.
In order to rescale the different graphs of $\kappa$ as a function of $\lambda$. 
Admitting that the GEV statistics will be correctly sampled if the bin length $m$ scales as the exit time necessary, we rewrite (\ref{conj}) as
$$
\epsilon^2\log(m) \propto \Delta V\,.
$$
Figure~\ref{langevin} displays the results for the minima of observable $X$ for several combinations of bin length $m$ and noise intensity $\epsilon$, while varying the quality of the statistics by changing the number $n$ of bins considered.
It clearly shows that the location of the zero of $\kappa$ is statistically well defined through the proposed rescaling, since fixing $\lambda$, i.e. $\Delta V$ comes to specifying the physical conditions.
The presence of the logarithm for $m$ and the square for $\epsilon$ in the scaling relation points out intrinsic limitations to the identification of $\lambda_{\rm c}$ through the condition $\kappa=0$, while showing that good estimates can be obtained with limited data.
 
\section{Final Remarks and Perspectives\label{S-Conc}}

In this paper we have introduced a new method for detecting critical transitions based on the properties of extreme fluctuations in recorded time series of relevant observables measured on the system of interest.
We have shown that the change in the  shape parameter of the Generalized Extreme Value distribution is a valid criterion for locating the bifurcation while the probability that the system experiences a transition is still low in view of the available length of the data sample.
The main advantage of this method when compared to others based on the property of the bulk statistics is to provide a value for the critical parameter at which the bifurcation is likely to occur.
In contrast, as we have seen through the examples provided, variance and skewness can be used as qualitative warning but there is no way of using them for a quantitative estimation of critical parameters.

We have tested the capabilities of the method on transitional plane Couette flow, the non trivial bifurcation structure of which has been well established by previous studies.
We have proposed that the value of shape parameter $\kappa=0$ be the signature of the transition and to locate the transition at the corresponding value control parameter, hence the global stability threshold $R_{\rm g}$.
The method has indeed been able to provide estimates for $R_{\rm g}$, most significantly at low resolution with the largest possible amount of data, but also at medium and high resolution, with shorter time series. 

The statistical analysis subsequently carried out using toy-models has confirmed that the methodology was able to provide warnings efficiently, in that $\kappa=0$ marks a transition where the probability of transition becomes significant.
More theoretical work is however necessary to further justify the proposed criterion, especially regarding the variation of this probability with the control parameter, say $\lambda$, so that the extrapolation of the function $\kappa(\lambda)$ to zero effectively gives a reliable critical value $\lambda_{\rm c}$.
In particular, the amount of data necessary to this estimate remains to be assessed, however the continuous increase of data collection capabilities in experiments and the ever growing computational power support some optimism in this matter.  

The approach followed in this exploratory paper offers a novel avenue to tackle the problem of detecting bifurcations in complex systems and is hoped to trigger the search for rigorous justifications from theoretical studies and further confirmation from other experiments.
As a matter of fact, the path from empirical results on plane Couette flow to related toy models can be reproduced in many other systems where multi-stability problems are of interest, vis. in climate sciences, the stability of the overturning meridional circulation \cite{broecker1997thermohaline,lucarini2012bistable} or the problem of the  transition to a snowball planet \cite{hoffman2002snowball}.
For this reason we believe that,  in the future, the method will be accessible to a vast public of different scientific disciplines.

\paragraph{Acknowledgments}
This research has received funding from the European Research Council under the European Community Seventh Framework Programme (FP7/2007-2013) / ERC Grant agreement No. 257106 and  the HPC-EUROPA2 project (project number: 228398) with the support of the European Commission Capacities Area -- Research Infrastructures Initiative. DF acknowledges  CNRS for a post-doctoral grant.

\bibliographystyle{model1-num-names}
\bibliography{couettebib}

\end{document}